\documentclass[11pt, letterpaper, ]{article}

\usepackage[utf8]{inputenc}
\usepackage{amssymb}
\usepackage{multicol}
\usepackage[margin=2.9cm]{geometry}
\usepackage{setspace}
\usepackage{graphicx}
\usepackage{url}
\usepackage{color}
\usepackage{hyperref}
\graphicspath{ {images/} }
\usepackage{tikz}
\usepackage{amsmath}

\begin{document}
	
	\title{\textbf{Weyl's Quantifiers}}
	
	\author{}
	
	\date{Iulian D. Toader\thanks{Institute Vienna Circle, University of Vienna, iulian.danut.toader@univie.ac.at}}
	
	\maketitle
	
	
	\bigskip \bigskip 
	
	\begin{quote} \textbf{Abstract}: I argue against the predominant view of Weyl’s interpretation of the logical signs. Drawing on his correctness-first account of mathematical knowledge, I point out that, according to him, quantified statements generate conditional obligations to act in ways that expand the repository of correct judgments. This clarifies Weyl's reasons for rejecting the law of excluded middle, which have nothing to do with what has been attributed to him by the predominant view. I also offer some preliminary thoughts on how to understand conditional obligations generated by statements with nested quantifiers. \end{quote}
	
	\bigskip \bigskip 
	
	\section{Introduction}
	
	In one of his last publications, David Charles McCarty confessed: 
	
	\begin{quote} 
		
	I have never taken the contributions of Hermann Weyl to the foundations of intuitionistic mathematics at all seriously. For one thing, he was a turncoat, the Benedict Arnold of mathematical intuitionism. ... Weyl embraced (a ragged selection from) Brouwer's ideas with his usual grandiose fanfare. However -- by 1924 or so -- he had slunk away from them ignominiously. ... As if he had not inflicted enough damage upon Brouwer's intuitionism, Weyl set out to clarify the murk of the Brouwerian philosophy by folding into it liberal doses of the intellectually poisonous mire of Husserlian phenomenology. In that, Weyl only compounded his treachery and, like his latterday followers [Van Atten, D. van Dalen, \& Tieszen 2002], ramped up the crazy. Without question, there is a hell for mathematicians -- and Hermann Weyl now lies deep in its Ninth Circle. (McCarty 2021, 316)
		
	\end{quote}
	
	\noindent There are good reasons for some of these claims. Unfortunately, the secondary literature remains largely ineffective in suppressing the intellectual toxins that McCarty, and many other philosophers, detected in Husserlian phenomenology. Most attempts to elucidate the effects of the latter on Weyl's philosophical thinking have been just as successful.\footnote{In all fairness, Weyl is a hard nut to crack. For some of my own attempts, see Toader 2013, 2014, 2021.} But as far as I can see, some other of McCarty's claims are unwarranted. For one thing, it's not clear that what he saw as treachery can really be held against Weyl. It seems to me that, like his character Aristides, Weyl treasured honesty more than loyalty.\footnote{Surely some kinds of treason are morally justified, even obligatory (see e.g. Fabre 2020).} But I will pass no judgment on that, though I find McCarty's vision of the mathematical underworld quite amusing, just as he may have intended it.
	
	In any case, McCarty followed up his confession with more specific criticism: 
	
	\begin{quote} 
		
	Weyl put forward his own, idiosyncratic interpretations of the logical signs. ... Weyl took it to follow from these strange and conflicting reinterpretations of quantifiers and negation that the law of the excluded middle is invalid. ... Weyl has in no way here described any sort of counterexample to the law of the excluded middle or any sort of good reason for thinking it invalid. That there are `instances' of TND that are meaningless because the substituends were meaningless to start with, e.g., (with apologies to Lewis Carroll) `Either the slithy tove did gyre or the slithy tove did not gyre' shows nothing about the validity or invalidity of the law. One can take any standard expression of a valid logical scheme and substitute uniformly into its sentential variables nonsensical strings of words and get further nonsense out! Surprise!'' (op. cit., 319)
		
	\end{quote}
	
	\noindent This kind of overhasty criticism calls for a careful reconsideration of Weyl's interpretation of the logical signs, which I want to undertake in this paper. As we will see in section 2, the predominant view, which I will call \textit{the Ramseyan view} for it seems to have been first proposed by Ramsey, conflates Weyl's interpretation with that defended by Hilbert. Followed by everyone, including McCarty, the Ramseyan view effectively propagated a \textit{mis}interpretation, and allowed commentators to speak of ``the particular way in which Hilbert and Weyl demarcated meaningful and meaningless'' statements (Raatikainen 2003, 173), as if they both demarcated them in one and the same way. In fact, while he maintained that transfinite axioms are not proper judgments because they do not represent any state of affairs, Weyl never thought that they were meaningless, although he did point out (rightly or wrongly) that Hilbert himself considered them to be meaningless. For Weyl, as I will reconstruct his interpretation, the meaning of quantified statements cannot be given by their representation of any state of affairs -- there are no states of affairs that they can represent -- but only by their normative function. 
	
	Some likely historical sources for Weyl's interpretation are Schlick's 1918 \textit{Allgemeine Erkenntnislehre} and Alexander Pfänder's 1919 \textit{Logik} (Sundholm 1994, 125sq), or Schopenhauer's 1859 \textit{Die Welt als Wille und Vortstellung} (Detlefsen 2011, 90). It is entirely reasonable to believe that Weyl was fully aware of these works. It seems no less reasonable to assume that Weyl was also acquainted with Husserl's phenomenology of values, and in particular his lectures on axiology (Husserl 1988), delivered in the winter semester 1908/09 and the summer semesters 1911 and 1914 in Göttingen, where Weyl was based until he took up his position at the ETH in Zurich, in 1913. Despite McCarty's reservations, I think it would be very important and extremely interesting to investigate Husserlian axiology as a basis for Weyl's interpretation of the logical signs. But this historical investigation will have to be pursued elsewhere. Here I want to focus on Weyl's interpretation, quite independently of any sources that might have helped him to develop it.\footnote{For discussion of Husserlian axiology, see Mulligan 2004. Weyl was certainly familiar with the \textit{Logische Untersuchungen}, where Husserl gives a preliminary account of values (see Mulligan 2022), but it is an open question if that would have been sufficient as a basis for Weyl's interpretation of the logical signs.}
	
	More specifically, I will argue, in section 3, that what distinguishes Weyl's quantifiers from the usual logical ones can be explained by closely considering the normative function that I believe he assigned to quantified statements. He described this function in his usual (or, rather, unusual) metaphorical style, which has given much trouble to philosophers, from logical empiricists like Hans Hahn (see Hahn 1928) to mathematical intuitionists like McCarty. But if one looks harder, one can see through the metaphors, quite clearly, that the rules of Weyl's quantifiers are not the familiar rules of the usual logical quantifiers. Once this is well understood, the invalidity of the law of excluded middle in transfinite mathematics follows then, exactly as Weyl maintained, i.e., as an immediate consequence of his interpretation of the logical signs, rather than because the law would allow meaningless instances as counterexamples, as McCarty erroneously thought Weyl claimed.
	
	Taken seriously, and apart from the grandiose rhetoric of the 1920s, Weyl's interpretation implies that transfinite mathematics is a normative system, in the sense that its quantified statements are values that generate norms for mathematical action. These norms, as we will see, are understood deontically, as obligations to act in ways that expand the repository of what Weyl thought was the only source of mathematical knowledge: correct judgments, or judgments that correctly represent states of affairs. But can one take Weyl's interpretation seriously?\footnote{Also, \textit{should} one consider Weyl's view apart from his rhetoric, the knack for which he undoubtedly picked up from Hilbert (as noted in Martin-Löf 2008, 246)? For discussion of rhetoric in mathematical practice and the foundations of science, see e.g. Reyes 2004, Kragh 2015.} To make sense of his conception of transfinite mathematics as a normative system, one must fully explain the normative function of quantified statements. One particularly difficult aspect that needs to be better understood is the structure of obligations generated by statements with nested quantifiers. I will outline my preliminary thoughts on this in section 4, before concluding the paper.

	\section{Against the Ramseyan view}
	
	How did Weyl understand transfinite mathematics? To answer this question, I will present a new reconstruction of his interpretation of the logical signs, which he introduced in 1921, an interpretation that there is no reason to think he ever abandoned.\footnote{Indeed, Weyl kept referring back to it. See, for instance, Weyl 1929, 1940.} But this reconstruction is not what has been typically attributed to him. Consider the following analysis:
	
	\begin{quote} 
		
	Weyl's reasons for rejecting excluded middle were quite different from Brouwer's. Weyl seems to have accepted it for propositional calculus. His reason for rejecting it for the predicate calculus ... was that quantified sentences, not being meaningful, cannot be meaningfully negated; and hence that no instance of excluded middle is meaningful. Now it is not even clear that this should be regarded as a rejection of excluded middle, which is surely only meant to apply to meaningful sentences. We do not count as rejecting it if we merely refuse to accept `Either og ur blig or not og ur blig' on the grounds that `Og ur blig' is meaningless. In contrast [to intuitionism], the kind of approach favored by Weyl in 1920, and later by Ramsey, is very much in line with formalist thinking. Quantified sentences are thought of as devices, strictly meaningless in themselves, that allow for the manipulation of meaningful sentences. (Price 2011, 155)
		
	\end{quote}
	
	\noindent Just like McCarty's criticism (quoted above), Huw Price's analysis is mistaken. Quantified statements, for Weyl, are not meaningless, so he could not have rejected the law of excluded middle on account of their meaninglessness. Further, consider also Wilfried Sieg's contention that ``Weyl’s views are, in some important respects (the understanding of quantifiers is one such point) close to the finitist standpoint.'' (Sieg 2013, 117) This is inaccurate as well, for Weyl's understanding of the quantifiers was actually far enough from Hilbert's finitist standpoint. This widespread misinterpretation originated, as far as I have been able to determine, in Ramsey's paper, ``Mathematical Logic'', so I call it \textit{the Ramseyan view}:
	
	\begin{quote} 
		
	[A]nother more fundamental reason is put forward for denying the Law of Excluded Middle. This is that general and existential propositions are not really propositions at all. ... The foundation on which this rests [is] the view that existential and general propositions are not genuine judgments... \textit{Hilbert shares Weyl's opinion that general and existential propositions are meaningless... We must begin with what appears to be the crucial question, the meaning of general and existential propositions, about which Hilbert and Weyl take substantially the same view.} Weyl says that an existential proposition is not a judgment, but an abstract of a judgment, and that a general proposition is a sort of cheque which can be cashed for a real judgment when an instance of it occurs. Hilbert, less metaphorically, says that they are ideal propositions, and fulfill the same function in logic as ideal elements in various branches of mathematics. (Ramsey 1926, my emphasis) 
		
	\end{quote}
	
	\noindent Ramsey was right: the meaning of quantified statements \textit{is} the crucial question. Indeed, Weyl maintained that such statements are not genuine or real judgments, in the sense that they do not represent states of affairs, and Hilbert considered them as ideal statements. But whereas Ramsey reported (rightly or wrongly) that Hilbert took it to be the case that existential and universal statements are as such meaningless, he (Ramsey) failed to see that Weyl never proposed or endorsed this view: although they are not what he called proper judgments, quantified statements are nevertheless meaningful, in a sense to be presently clarified. 
	
	To some extent, Hilbert may have been influenced by Weyl's interpretation of quantified statements. Dirk van Dalen, among others, claimed that this was indeed the case: ``Weyl's conception of existential and general statements has influenced the treatment of the finitary viewpoint in arithmetic.'' (van Dalen 1995, 158) He also contended that Hilbert took such statements to have ``a hypothetical finitary meaning'' (\textit{loc. cit.}), or as Hilbert actually put it, they ``can be interpreted finitarily only in a hypothetical sense'' (Hilbert and Bernays 1934, 32). However, despite what Ramsey (and Weyl) reported, it is undeniable that Hilbert recognized the need to show them meaningful, to show that an existential statement, for example, is meaningful ``independently of its [finitary] interpretation as a partial judgment'' (\textit{op. cit.}, 37). But as far as I can see, there is no evidence that Hilbert came close to Weyl's interpretation of quantification as reconstructed in this paper.
	
	Furthermore, Weyl never maintained that quantified statements fulfill the same function in logic as ideal elements do in mathematics. Hilbert described his own view as follows: ``In my proof theory, the transfinite axioms and formulae are adjoined to the finite axioms, just as in the theory of complex variables the imaginary elements are adjoined to the real, and just as in geometry the ideal constructions are adjoined to the actual.'' (Hilbert 1923, 1144) The function of transfinite axioms is to increase the simplicity and fruitfulness of finitary mathematics, just like the introduction of imaginary numbers, say, was meant to increase the simplicity and fruitfulness of real analysis. But, again, Weyl never endorsed this view. Regarding the fundamental theorem of algebra, he thought that ``its use should be avoided as long as possible.'' (Weyl 1932, 46). For him, only proper judgments could fulfill a function in logic. Quantified statements fulfill their function not as judgments, but as a kind of values. And as we will see, their primary function is not, or perhaps not primarily, to increase simplicity and fruitfulness, but to generate obligations to act in certain ways.
	
	Once we realize this, we can begin to understand where the Ramseyan view went wrong, and we can see that refusing to accept ``either the slithy tove did gyre or the slithy tove did not gyre'' on the grounds that ``the slithy tove did gyre'' is meaningless, just like refusing to accept ``either og ur blig or not og ur blig'' on the grounds that ``og ur blig'' is meaningless, has nothing to do with Weyl's actual reasons for rejecting the validity of the law of excluded middle.\footnote{It is rather surprising that one could think such examples enough to doubt the cogency of Weyl's reasons. McCarty's instance of the excluded middle is not even meaningless. Hence, I take it, his apologies to Carroll. For as Humpty Dumpty explains to Alice, ``slithy'' means ``lithe'' and ``slimy'', or smooth and active, and ``toves'' are like badgers or lizards, with long hind legs, and short horns like a stag, that lived chiefly on cheese. ``Gyre'' is derived from ``gyaour'' and means ``to scratch like a dog.'' For details, see Gardner 2000.} 
	
	The refutation of the Ramseyan view starts by clarifying a central element of Weyl's philosophy of mathematics: the notion of proper judgment, or what he called judgment in a proper sense. As early as 1918, in his \textit{Das Kontinuum}, he explained this notion as part of his characterization of what he called elementary inferences in mathematics: a proper judgment is a judgment that represents a state of affairs, and elementary inferences are inferences that involve only such judgments (Weyl 1918, 17). Weyl understood inferential correctness to depend essentially on representational correctness, and defined entailment in terms of an order relation between correctly represented states of affairs. For example, to borrow an illustrious example from Locke, reasoning from ``AB is an arch of a circle'' to ``AB is less than the whole circle'' is inferentially correct in virtue of the order relation, and more exactly the containment relation, between these two correctly represented states of affairs. According to Weyl, a genuine mathematical proof should deploy only inferentially correct reasoning and, thus, should contain only representationally correct judgments. He explicitly endorsed a correctness-first account of knowledge: ``correctness ... remains throughout the ultimate source of knowledge.'' (\textit{op. cit.}, 11) As the term is used here, ``correctness'' designates an epistemic aspect of the representational capacity of judgments whose semantic aspect is truth. What led him to this account was his early acquiescence to a phenomenological epistemology, which conceived of correctness as a norm, in a deontic sense: correct judgments are judgments that we ought to judge. For Weyl then, a genuine proof requires that we fully discharge our unconditional obligation to judge correctly.\footnote{On this epistemological background, Weyl criticized Dedekind's conception of proof, which he (Weyl) thought promoted a notion of inferential correctness that is independent of representational correctness and, therefore, is epistemologically problematic. For some discussion of this criticism, see Toader 2016. The obligation to judge correctly may be regarded as an epistemic obligation, although this notion typically refers to obligations in matters of belief or opinion, rather than judgment. For discussions of the notion of correctness, as understood in the phenomenological tradition, see Mulligan 2017, and Textor 2019, 2022.} 
	
	This account of proof has important consequences for what Ramsey perceived as the crucial question of the meaning of quantified statements. Since they do not represent states of affairs, quantified statements are not proper judgments, so they cannot be part of what Weyl would regard as a genuine proof. Still, he did not think that quantified statements have no contribution to knowledge. Indeed, it is precisely the function they have in transfinite mathematics that supports and explains their contribution. Weyl's description of this function is given in his \textit{erlösende Wort}, in two fragments that are among the most often quoted texts from his entire work. Here is the first fragment:
	
	\begin{quote} 
		
	An existential statement -- something like ``there is an even number'' -- is not a judgment in a proper sense, one that asserts a state of affairs. Existential states of affairs are an empty invention of logicians. ``2 is an even number'': this is a real judgment that gives expression to a state of affairs; ``there is an even number'' is only a judgment-abstract obtained from this judgment. Taking knowledge to be a valuable treasure, the judgment-abstract is a paper that indicates the presence of a treasure without disclosing where it is. Its only value can lie in its ability to get me to search for the treasure. The paper is worthless so long as a real judgment that stands behind it, such as ``2 is an even number'', is not provided. (Weyl 1921, 54sq, translation amended)
		
	\end{quote}
	
	\noindent Against the Ramseyan view, let me emphasize that Weyl did not say that existential statements are meaningless. They would be meaningless if one mistook them for proper judgments. But he explicitly maintained that they cannot be taken as proper judgments, for they do not represent any state of affairs. What relations do existential statements, as judgment-abstracts, have to our obligation to judge correctly? How can they serve the provision of genuine proofs? Weyl wrote that existential statements are \textit{indications} of the existence of proper judgments. As such, he considered them valuable, insofar as they prompt us to act in certain ways -- to search for the proper judgments whose existence they indicate. Given Weyl's emphasis on our obligation to judge correctly, I think it is natural to understand indications deontically as generating a further obligation: we ought to search for correct judgments. This obligation is conditional on existential quantification. I will revisit this point in the next section, where I give a more precise reconstruction of Weyl's reasoning in the first fragment, which will clearly show that the rules of Weyl's existential quantifier are different than the usual logical rules.
	
	A defender of the Ramseyan view might insist that valuable indications could still be meaningless. But honestly, I fail to see how meaningless indications could prompt one to act in the way Weyl described, or any way, for that matter. My view is that it is exactly the normative function that he attributed to existential statements that makes them meaningful. This is where I think the Ramseyan view went wrong in this case.
	
	Let's turn now to the second fragment of Weyl's \textit{erlösende Wort}:
	
	\begin{quote} 
		
	Just as little is the general statement ``each number has property $\mathrm{F}$'' -- e.g., ``for any number \textit{m}, \textit{m} + 1 = 1 + \textit{m}'' -- a real judgment, but rather a general instruction for [recovering such] judgments. Based on this instruction, if I come across an individual number, e.g., 17, I can redeem a real judgment, that is: 17 + 1 = 1 + 17. To use another image: comparing knowledge to a piece of fruit and the intuitive realization of knowledge to its consumption, then a general statement is to be compared to a hard shell filled with fruit. It is, obviously, of some value, however, not as a shell by itself, but only for its fruit content. It is of no use to me as long as I do not open it and actually take out the piece of fruit and eat it. (\textit{loc. cit.})
		
	\end{quote}
	
	\noindent Weyl did not say that universal statements are meaningless, either. Just as above, one could consider them meaningless if one mistook them for proper judgments. But he again emphasized that they are not to be taken as proper judgments. What roles do universal statements play in mathematics then, on Weyl's correctness-first account of mathematical knowledge? How can \textit{they} serve the provision of genuine proof? Weyl wrote that universal statements are \textit{instructions} for the recovery of proper judgments. As such, he considered them valuable insofar as they prompt us to act in certain ways -- to recover proper judgments. Again, given Weyl's emphasis on our obligation to judge correctly, I think it is also natural to understand instructions normatively as generating a further obligation: we ought to recover correct judgments. This obligation is conditional on universal quantification. I will revisit this point as well in the next section, where I give a more exact reconstruction of Weyl's reasoning in the second fragment, which will clearly show that the rules of Weyl's universal quantifier are different than the usual logical rules.
	
	Neverthless, the defender of the Ramseyan view might still want to insist that instructions could be valuable, but meaningless. Again, I just cannot see how meaningless instructions could prompt one to act in the way Weyl described, or any way, for that matter. Once more, my view is that it is exactly the normative function that Weyl attributed to universal statements that makes them meaningful. This is where the Ramseyan view went wrong in this other case as well. 
	
	I turn now to explaining in more detail why I think that, according to Weyl, quantified statements generate obligations, and how this clarifies his interpretation of the logical signs and his actual reasons for rejecting the law of excluded middle. Of course, it's difficult to say that the reconstruction I propose is entirely faithful to Weyl's thought, but I think it comes close to it, and in any case, it is preferable to the Ramseyan view.

	\section{Transfinite mathematics as a normative system}
	
	There may be different ways of understanding Weyl's apparently strange and metaphorical remarks in the two fragments quoted in the previous section. But I think that they are well understood if presented as follows:
	
	\begin{quote} 
		
	1. We value proper judgments.
		
	2. If we assert an existential statement, we give an indication of the existence of a proper judgment.
		
	3. Thus, we value an existential statement because it prompts us to search for a proper judgment.
		
	4. If we assert a universal statement, we give an instruction for recovering a proper judgment whenever an object of the domain is given.
		
	5. Thus, we value a universal statement because it prompts us to recover a proper judgment whenever an object of the domain is given.
		
	\end{quote}
	
	\noindent This reconstruction explains why Weyl believed that quantified statements are valuable even though they do not represent any state of affairs. It also explains how the value of these statements derives from that of proper judgments. Essential in this derivation of their value is the ability of quantified statements to prompt us to act in certain ways. If they did not have the ability to prompt us to search for proper judgments, then existential statements would have no value. If they did not have the ability to prompt us to recover proper judgments, then universal statements would have no value, either. 
	
	However, this reconstruction does not take into account, in claim 1, Weyl's correctness-first epistemology, and thereby fails to clarify the nature of the prompting ability of quantified statements. As a consequence, the derived value he took quantified statements to have is not adequately explained in claims 3 and 5. As we have seen, he believed that the view according to which correctness is the only source of knowledge implies that genuine proofs require that we discharge our obligation to judge correctly. I take this to suggest that the ability of quantified statements to prompt us to act in certain ways is best conceived of normatively, as an ability to generate an obligation to act in those ways. Without discharging this conditional obligation, i.e., without searching for a \textit{correct} judgment when we assert an existential generalization, or without recovering a \textit{correct} judgment when we assert a universal generalization, we would not be in a position to disclose or redeem correct judgments, and thus we would fail to discharge our unconditional obligation to judge correctly. Taking this into account modifies the above reconstruction as follows:
	
	\begin{quote} 
		
	1'. We ought to judge correctly because we value correct judgments.
		
	2'. If we assert an existential statement, we ought to search for a correct judgment.
		
	3'. Thus, we value an existential statement because it generates an obligation to search for a correct judgment. 
		
	4'. If we assert a universal statement, we ought to recover a correct judgment whenever an object of the domain is given.
		
	5'. Thus, we value a universal statement because it generates an obligation to recover a correct judgment whenever an object of the domain is given. 
		
	\end{quote}
	
	\noindent But there is still a problem with claim 4'. Recovering a proper judgment whenever an object of the domain is given is not always possible, because we cannot run through an infinite set (Weyl 1921, 54). To address this point, Weyl proposed that a universal statement is to be regarded as a special kind of existential statement, one that asserts the existence of an essence. Thus, an instruction for recovering a correct judgment should be reconceived as an indication of the existence of an essence. This suggests the replacement of claims 4' and 5' by the following ones:
	
	\begin{quote} 
		
	4''. If we assert a universal statement, we ought to search for an essence in order to recover a correct judgment.
		
	5''. Thus, we value a universal statement because it generates an obligation to search for an essence in order to recover a correct judgment.
		
	\end{quote}
	
	\noindent This reconstruction explains not only why Weyl believed that quantified statements are valuable even though they do not represent any state of affairs, and how their value derives from that of correct judgments, but it clarifies as well that what is absolutely essential in this derivation is the ability of quantified statements to generate an obligation to act in certain ways. If they did not have an ability to generate an obligation to search for correct judgments, existential statements would have no value. If they did not have an ability to generate an obligation to recover correct judgments, universal statements would have no value, either. Their ability to generate such obligations underlies their normative function. This function makes them meaningful, and supports and explains their contribution to knowledge: they serve what Weyl saw as a fundamental epistemological task of the mathematician, i.e., to discharge the unconditional obligation to judge correctly.
	
	This understanding of the function of quantified statements implies that Weyl's quantifiers have different rules than the familiar rules of the usual logical quantifiers. These rules are expressed in claims 2' and 4''. They may be regarded as deontic rules for quantifier elimination, rules that prescribe how we ought to act when we assert quantified statements. Reconstructed in this way, Weyl's view is neither strange, nor really metaphorical. Unlike Hilbert, who understood quantified statements as partial assertions or ideal statements that may have a function in logic, Weyl assigned them the normative function just described. But if this reconstruction is accurate, as I believe it to be, then it follows that he never conceived of transfinite mathematics as a game with meaningless symbols, despite the fact that he did attribute this view (rightly or wrongly) to Hilbert. Rather, Weyl saw it as a normative extension of elementary mathematics, as a system of conditional obligations.\footnote{Conditional obligations have been studied, of course, since at least the 1950s (see e.g. von Wright 1956), and various problems and puzzles with both unconditional and conditional obligations, especially contrary-to-duty obligations, have been extensively discussed (see e.g. Chisholm 1963, van Fraassen 1972, and most recently Fine 2024). My reconstruction would be strengthened by solutions to at least some such problems as might arise for Weyl's conception of transfinite mathematics, but this will have to be addressed in future work.}
	
	On this conception of transfinite mathematics, one can now understand Weyl's reasons for maintaining that the law of excluded middle is invalid. The Ramseyan view, recall, had it that ``quantified sentences, not being meaningful, cannot be meaningfully negated; and hence ... no instance of excluded middle is meaningful.'' On this view, excluded middle fails because negated quantified statements are meaningless, and this is because quantified statements themselves are meaningless. However, as I have argued, this is not Weyl's view. For he simply did not take quantified statements to be meaningless, although he did explicitly maintain that it would be utterly meaningless to negate such statements:
	
	\begin{quote} 
		
	Our theory of general and existential statements is not intuitively vague, something that, amongst  other things, becomes apparent from the fact that it leads immediately to important, rigorously clear consequences. Above all, that it is completely meaningless to negate this sort of statements. Thus, the possibility to formulate an ``axiom of the excluded middle'' for them disappears. (Weyl 1921, 56, translation amended)
		
	\end{quote}
	
	\noindent In 1929, in his first lecture at the Rice Institute, Weyl repeated that it is ``evidently impossible and without meaning'' to negate quantified statements (Weyl 1929, 247). What is the best way to understand this claim as an immediate consequence of his understanding of quantified statements as having a normative function? My answer is that what Weyl meant is that, although quantified statements are meaningful, it would be impossible and meaningless to \textit{logically} negate them. For he thought that only proper judgments could be subject to logical negation, but quantified statements should not be mistaken for proper judgments. According to Weyl then, the law of excluded middle fails in transfinite mathematics because logically negated quantified statements are meaningless. Again, this is not because the quantified statements themselves are meaningless. Rather, it is because, as valuable instructions and indications that generate obligations to act in certain ways, they cannot be logically negated. Negation, obeying the familiar logical rules, does not apply to this kind of meaningful statements. 
	
	In the same lecture to the Rice institute, Weyl further noted: ``With regard to ... the usage of the terms `all' and `any,' I think one does not hit quite the right spot by referring to the validity or invalidity of the principle of the excluded middle." (\textit{loc. cit}). The right spot to be hit, ``the point on which the matter hinges'', Weyl specified, is ``the fact that the negation cannot be carried out'' for quantified statements. If my reconstruction above is accurate, then Weyl did not take logical revision, i.e., the rejection of the law of excluded middle in transfinite mathematics, to tell us anything about the meaning of quantifiers and negation. Rather, it was the other way around: semantic revision, i.e., his own interpretation of quantifiers and negation, justifies the logical revision. At the same time, it should be clear enough by now that Weyl did not provide a counterexample to the law of excluded middle by uniform substitution into its sentential variables of nonsensical strings of words like ``og ur blig'' or, perhaps, ``the slithy tove did gyre''. Rather, he thought that there were meaningful statements that would do the job. And they do the job if interpreted as my reconstruction suggests they should.

	\section{Nested quantifiers}
	
	As a normative system, Weyl's conception of transfinite mathematics may, of course, admit its own kind of negation -- different than the usual logical negation -- one that could be meaningfully applied to quantified statements without validating the law of exclude middle. Moreover, his restriction of the applicability of logical negation to proper judgments seems to further imply that transfinite mathematics cannot be \textit{logically} inconsistent, though it could perhaps be \textit{deontically} inconsistent. I discuss these issues elsewhere, in the context of Weyl's argument, against Hilbert, for the dispensability of consistency proofs. In this section, I want to focus on another important question: what is the structure of obligations that are conditional on quantification? More specifically, how should we conceive of the obligations generated by statements with nested quantifiers? What sense can we make of Weyl's remarks about statements whose quantifiers occur within the scope of other quantifiers, if my reconstruction of his interpretation is adopted instead of the Ramseyan view? Here are some preliminary thoughts on this matter.
	
	When turning to nested quantifiers, Weyl first mentioned the possibility to existentially generalize a universal statement: ``It is possible to draw an abstract not only from a [proper] judgment, but also from an instruction for [recovering a proper] judgment.'' (Weyl 1921, 57, translation amended) To explain how to draw an abstract from an instruction, he gave the following example: ``every number \textit{m} stands in relation \textit{R} to 5''. The abstract, in this case, is ``There is a number \textit{n} such that every number \textit{m} stands to it in the relation \textit{R(m, n)}.'' Weyl considered this abstract unproblematic, and understood it as the indication of an instruction for recovering a proper judgment. Taking again into account his correctness-first epistemology, my reconstruction explains what  obligations are generated by existential generalizations of universal statements and why we value the latter: 
	
	\begin{quote}  
		
	6. If we assert an existential generalization of a universal statement, we ought to search for an instruction for recovering a correct judgment, and then we ought to recover a correct judgment whenever an object of the domain is given.
		
	7. Thus, we value an existential generalization of a universal statement because it generates an obligation to search for an instruction for recovering a correct judgment, an instruction which in turn generates an obligation to recover a correct judgment whenever an object of the domain is given.
		
	\end{quote}
	
	\noindent Secondly, Weyl turned to universal generalizations of existential statements, and explained what he saw as a problem:
	
	\begin{quote} 
		
	By contrast, an instruction for [recovering] judgment abstracts is simply nothing (``das reine Nichts''), so long as it is not backed by an instruction for [recovering] the real judgments from which it [i.e., the instruction for recovering judgment abstracts] has been obtained as an abstract. Example: For every number \textit{m} there is a number \textit{n} such that the relation \textit{R(m, n)} holds between them. We must in truth be dealing here with an abstract from a judgment instruction. Which judgment instruction? (\textit{loc. cit.})   
		
	\end{quote}
	
	\noindent Before we see how Weyl answers this question, we have to understand why he asks \textit{this} question. Why must we be dealing here with an abstract from a judgment instruction, rather than with an instruction for recovering a judgment abstract? Key to this is his claim that universal generalizations of existential statements, like ``For every number \textit{m} there is a number \textit{n} such that the relation \textit{R(m, n)} holds between them.'' must be interpreted as ``There is a \textit{law} $\phi$ such that for every number \textit{m} the relation $R$ holds between $m$ and $\phi(m)$.'' More generally, any instruction for recovering a judgment abstract must be interpreted as the indication of a law that backs the instruction for recovering a proper judgment. This reduces the case of universal generalizations of existential statements to the unproblematic case of existential generalizations of universal statements.  Taking once more into account Weyl's correctness-first epistemology, my reconstruction explains what obligations are generated by universal generalizations of existential statements and why we value the latter:
	
	\begin{quote} 
		
	8. If we assert a universal generalization of an existential statement, then we ought to search for a law in order to recover a correct judgment.
		
	9. Thus, we value a universal generalization of an existential statement because it generates an obligation to search for a law in order to recover a correct judgment.
		
	\end{quote}
	
	\noindent Weyl's reduction explicitly assumes that an instruction for recovering judgment abstracts can only be obtained as an abstract from an instruction for recovering correct judgments. But why should this be so? In other words, why do universal generalizations have to be interpreted in terms of the existence of a law or, as we have seen in the previous section, the existence of an essence? Did Weyl offer any justification for this assumption (expressed by claim 8)? It is surely impossible for us to run through infinite sets, and maybe as he contended this would be impossible even for God, but why exactly does this make positing essences or laws indispensable for interpreting and understanding universal generalizations? On my reconstruction, it is clear why positing essences or laws is indispensable for discharging the obligations generated by universal generalizations. For suppose we reject claim 8, and so we don't take an instruction for recovering a judgment abstract as an indication of a law that backs an instruction for recovering a correct judgment. Then the obligations generated by a statement with an existential quantifier within the scope of a universal quantifier would be different: 
	
	\begin{quote} 
		
	8'. If we assert a universal generalization of an existential statement, we ought to recover an indication of the existence of a correct judgment, and then we ought to search for a correct judgment.
		
	9'. Thus, we value a universal generalization of an existential statement because it generates an obligation to recover an indication of the existence of a correct judgment, an indication which in turn generates an obligation to search for a correct judgment.
		
	\end{quote}
	
	\noindent However, for precisely the epistemic reasons indicated by Weyl, it is doubtful that we can always discharge the obligation to recover an indication of the existence of a correct judgment. But then, without an indication of the existence of a correct judgment, the conditional obligation to search for a correct judgment would not be generated, and therefore we would at least sometimes fail to discharge our unconditional obligation to judge correctly. 
	
	\section{Conclusion}
	
	The Ramseyan view, which has been the predominant view of Weyl's interpretation of the logical signs, conflates the latter with Hilbert's. It maintains that they both had the same understanding of quantification, and they both demarcated meaningful from meaningless statements in one and the same way. In this paper, I have offered some reasons to resist this view. I have pointed out that a central element of Weyl's interpretation is his correctness-first account of mathematical knowledge. This account, according to which we are under an unconditional obligation to judge correctly, suggests an understanding of quantified statements as generating further, conditional obligations to act in ways that expand the repository of correct judgments. I have argued that this justifies the claim that Weyl conceived of transfinite mathematics as a normative extension of elementary mathematics. This clarifies Weyl's semantic reasons for rejecting the validity of the law of excluded middle, which have nothing to do with the reasons attributed to him on the Ramseyan view. 
	
	My reconstruction of Weyl's interpretation is articulated by claims 1', 2', 3', 4'', and 5'' above. I have also offered some preliminary thoughts on how to understand his remarks on nested quantification, and the obligations generated by statements with quantifiers within the range of other quantifiers, as articulated by claims 6, 7, 8, and 9. I think that the reconstruction can further be fruitfully used in an attempt to clarify Weyl's remarks on choice sequences, and in particular, the fact, often noted, that for him, ``the universal quantifier ranges over lawless, the existential over lawlike sequences'' (van Atten, van Dalen, and Tieszen 2002, 220sq). If a universal statement cannot generate an obligation to recover a correct judgment, but only an obligation to search for a law in order to recover a correct judgment, then it follows, in particular, that universal quantification over lawless sequences cannot generate a obligation to recover a correct judgment about lawless sequences, but only an obligation to search for a law. What could be then recovered on the basis of a law is, of course, only a correct judgment about lawlike sequences. Note that this does not entail that there cannot be correct judgments about lawless sequences, but only that the universal quantifier does not generate any obligation to recover such judgments. More needs to be said, however, to provide a fully satisfactory account.

	\section*{Acknowledgments}

	Thanks to the editors for inviting me to contribute to this volume, to the referees for their comments, and to Kevin Mulligan for his thoughts on a previous lengthier draft of this paper. I am also indebted to the late David Charles McCarty for his encouragement and learnedness, which he always shared most kindly.
	
	\section*{References} 
	
	\noindent Chisholm, R. M. (1963) ``Contrary-to-Duty Imperatives and Deontic Logic'', in \textit{Analysis}, 24, 
	
	33-36
	
	\smallskip
	
	\noindent Detlefsen, M. (2011) ``Discovery, invention and realism: Gödel and others on the reality of
	
	concepts'', in \textit{Meaning in Mathematics}, ed. by J. Polkinghorne, OUP, 73-94
	
	\smallskip
	
	\noindent Fabre, C. (2020) ``The Morality of Treason'', in \textit{Law and Philosophy}, 39, 427-461
	
	\smallskip
	
	\noindent Fine, K. (2024) ``Conditional and Unconditional Obligation'', in \textit{Mind}, 133, 377-399
	
	\smallskip 
	
	\noindent Gardner, M. (2000) \textit{The Annotated Alice}, W. W. Norton, 1st ed. 1960
	
	\smallskip
	
	\noindent Hahn, H. (1928) ``Philosophie der Mathematik und Naturwissenschaft'', in
	\textit{Monatshefte der}
	
	\textit{Mathematik und Physik}, 35, 51-55
	
	\smallskip
	
	\noindent Hilbert, D. (1923) ``Die logischen Grundlagen der Mathematik'', in  \textit{Mathematische Annalen}, 
	
	88, 151–65. Eng. tr. in \textit{From Kant to Hilbert: A Source Book in the Foundations of} 
	
	\textit{Mathematics}, vol. 2, ed. by W. B. Ewald, OUP, 1996, 1134–1148
	
	\smallskip
	
	\noindent Hilbert, D. and P. Bernays (1934) \textit{Grundlagen der Mathematik I}, Springer, 2nd ed. 1968
	
	\smallskip
	
	\noindent Husserl, E. (1988) \textit{Vorlesungen über Ethik und Wertlehre} (1908–1914), ed. by U. Melle, 
	
	\textit{Husserliana} XXVIII, Kluwer
	
	\smallskip
	
	\noindent Kragh, H. (2015) ``Mathematics and Physics: The Idea of a Pre-Established Harmony'', in  
	
	\textit{Science \& Education}, 24, 515–527 
	
	\smallskip

	\noindent Martin-Löf, P. (2008) ``The Hilbert-Brouwer controversy resolved?'', in \textit{One Hundred Years} 
		
	\textit{of Intuitionism (1907-2007)}, ed. by M. van Atten, P. Boldini, M. Bourdeau, and G. 
	
	Heinzmann, Birkhäuser, 243-256

	\smallskip
	
	\noindent McCarty, D. C. (2021) ``Continuity in Intuitionism'', in \textit{The History of Continua:} 
	
	\textit{Philosophical and Mathematical Perspectives}, ed. by S. Shapiro and G. Hellman, OUP, 
	
	300-327
	
	\smallskip
	
	\noindent Mulligan, K. (2004) ``Husserl on the `Logics' of Valuing, Values and Norms'', in \textit{Fenomenologia}
	
	\textit{dalla Ragion Pratica}, ed. by B. Centi \& G. Gigliotti, Naples: Bibliopolis, 177-225
	
	\smallskip
	
	\noindent Mulligan, K. (2017) ``Brentano’s Knowledge, Austrian Verificationisms, and Epistemic 
	
	Accounts of Truth and Value'', in \textit{The Monist}, 100, 88-105
	
	\smallskip
	
	\noindent Mulligan, K. (2022) ``Logic, Logical Norms, and (Normative) Grounding'', in \textit{Bolzano's}
	
	\textit{Philosophy of Grounding}, ed. by S. Roski and B. Schnieder, OUP, 244-275
	
	\smallskip
	
	\noindent Price, H. (2011) \textit{Naturalism Without Mirrors}, ch. 7. Reprint of R. Holton and H. Price (2003) 
	
	``Ramsey on Saying and Whistling: A Discordant Note'', in \textit{Nous}, 37, 325-341
	
	\smallskip
	
	\noindent Raatikainen, P. (2003) ``Hilbert's Program Revisited'', in \textit{Synthese}, 137, 157-177
	
	\smallskip
	
	\noindent Ramsey, F. (1926) ``Mathematical Logic'', in \textit{The Foundations of Mathematics and other} 
	
	\textit{Logical Essays}, ed. by R. B. Braithwaite, 1931
	
	\smallskip
	
	\noindent Reyes, G. M. (2004) ``The Rhetoric in Mathematics: Newton, Leibniz, the Calculus, and the 
	
	Rhetorical Force of the	Infinitesimal'', in \textit{Quarterly Journal of Speech}, 90, 163–188
	
	\smallskip
	
	\noindent Sieg, W. (2013) \textit{Hilbert's Programs and Beyond}, OUP 
	
	\smallskip
	
	\noindent Sundholm, G. (1994) ``Existence, Proof and Truth-Making: A Perspective on the Intuitionistic
	
	Conception of Truth'', in \textit{Topoi} 13, 117-126
	
	\noindent Textor, M. (2019) ``Correctness First: Brentano on Judgment and Truth'', in \textit{The Act and}
	
	\textit{Object of Judgment Historical and Philosophical Perspectives}, ed. by B. Ball and C. 
	
	Schuringa, New York: Routledge, 129–150
	
	\smallskip
	
	\noindent Textor, M. (2022) ``That’s correct! Brentano on intuitive judgement'', in \textit{British Journal for}
	
	\textit{the History of Philosophy}, 31, 805-824
	
	\smallskip
	
	\noindent Toader, I. D. (2013) ``Concept formation and scientific objectivity: Weyl's turn against 
	
	Husserl'', in \textit{Hopos: The Journal of the International Society for the History of Philosophy} 
	
	\textit{of Science}, 3, 281-305
	
	\smallskip
	
	\noindent Toader, I. D. (2014) ``Why Did Weyl Think that Formalism's Victory Against Intuitionism 
	
	Entails a Defeat of Pure Phenomenology?'', in \textit{History and Philosophy of Logic}, 35,	198-208
	
	\smallskip
	
	\noindent Toader, I. D. (2016) ``Why Did Weyl Think that Dedekind's Norm of Belief in Mathematics 
	
	is Perverse?'', in \textit{Early Analytic Philosophy. New Perspectives on the Tradition}, ed. by S. 
	
	Costreie, Springer, The Western Ontario Series in Philosophy of Science, 80, 445-451
	
	\smallskip
	
	\noindent Toader, I. D. (2021) ``Why Did Weyl Think That Emmy Noether Made Algebra the Eldorado 
	
	 of Axiomatics?'', in \textit{Hopos: The Journal of the International Society for the History of}
	
	\textit{Philosophy of Science}, 11, 122-142
	
	\smallskip
	
	\noindent van Atten, M., D. van Dalen, and R. Tieszen (2002) ``Brouwer and Weyl: The Phenomenology 
	
	and Mathematics of the Intuitive Continuum'', in \textit{Philosophia Mathematica}, 10, 203–226
	
	\smallskip
	
	\noindent van Dalen, D. (1995) ``Hermann Weyl's intuitionistic mathematics'', in \textit{Bulletin of Symbolic}
	
	\textit{Logic}, 1, 145-169
	
	\smallskip
	
	\noindent van Fraassen, B. C. (1972) ``The Logic of Conditional Obligation'', in \textit{Journal of Philosophical}
	
	\textit{Logic}, 1, 417-438
	
	\smallskip
	
	\noindent von Wright, G. H. (1956) ``A note on deontic logic and derived obligation'', in \textit{Mind}, 65, 
	
	507-509 
	
	\smallskip
	
	\noindent Weyl, H. (1918) \textit{Das Kontinuum}. In \textit{Das Kontinuum und andere Monographien}, Chelsea  
	
	Publishing Company, 1960. Eng. tr. as \textit{The Continuum}, Th. Jefferson Univ. Press, 1987
	
	\smallskip
	
	\noindent Weyl, H. (1921) ``Über die Neue Grundlagenkrise der Mathematik'', in \textit{Mathematische Zeitschrift}, 
	
	10, 39–79. Eng. tr. in P. Mancosu (ed.) \textit{From Brouwer to Hilbert. The Debate on the}
	
	\textit{Foundations of Mathematics in the 1920s}, OUP, 86-118 
		
	\smallskip
	
	\noindent Weyl, H. (1929) ``Consistency in Mathematics'', in \textit{The Rice Institute Pamphlet}, 16, 245-265
	
	\smallskip
	
	\noindent Weyl, H. (1932) ``Topologie und abstrakte Algebra als zwei Wege mathematischen
	Verständnisses'' 
	
	in \textit{Unterrichtsblätter für Mathematik und Naturwissenschaften},  38, 177–88. Eng. tr. in 
	
	\textit{Levels of Infinity}, ed. by P. Pesic, New York: Dover, 33–48
	
	\smallskip
	
	\noindent Weyl, H. (1940) ``The Ghost of Modality'', in \textit{Philosophical Essays in Memory of Edmund}, 
	
	\textit{Husserl}, ed. by M. Farber, Cambridge University Press, 278-303

\end{document}